%% file: polyzetas.tex
\documentclass[a4paper, twoside, 10pt]{article}
\usepackage{amsmath,amssymb,amsthm}
\usepackage[francais,english]{babel} 
 
\def \AM{{\mathbb{A}}}

\def \CM{{\mathbb{C}}}

\def \NM{{\mathbb{N}}}

\def \QM{{\mathbb{Q}}}
\def \RM{{\mathbb{R}}}

\def \DC{{\mathcal D}}

\def \LC{{\mathcal L}}
\def \MC{{\mathcal M}}

\def \RC{{\mathcal R}}

\def\gal{\hbox{\rm Gal\,}}
\def\l{\lambda}

\def \BG{{\mathfrak B}}

\def \LG{{\mathfrak L}}

\def \aG{{\mathfrak a}}

\def \dG{{\mathfrak d}}
\def \eG{{\mathfrak e}}

\def \gG{{\mathfrak g}}

\def \iG{{\mathfrak i}}

\def \lG{{\mathfrak l}}
\def \mG{{\mathfrak m}}
\def \nG{{\mathfrak n}}

\def \rG{{\mathfrak r}}

\def \tG{{\mathfrak t}}
\def \uG{{\mathfrak u}}

\def \leq{\leqslant}
\def \geq{\geqslant}
\def\D{\Delta}
\def\wh{\widehat}
\def\spec#1{\hbox{\rm Spec}(#1)}
\def\ot{\mathop\otimes}
\def\eme{\raisebox{1ex}{\scriptsize {\`e}me}}
\def\emes{\raisebox{1ex}{\scriptsize {\`e}mes}}
\def\k{\Bbbk}

\newtheorem{prop}{Proposition}[section]

\newtheorem{defi}[prop]{D{\'e}finition}
\newtheorem{theo}{Th{\'e}or{\`e}me}

\newtheorem{conj}{Conjecture}

\usepackage{polyzetas} 
\evensidemargin \oddsidemargin 

\title{Torseurs associ{\'e}s {\`a} certaines relations alg{\'e}briques entre
  polyz{\^e}tas aux racines de l'unit{\'e}}
\author{Georges Racinet}
\date{version 3 --- 5 mars 2001\footnote{Les versions 1 et 2 sont
  respectivement du 4 et 11 d{\'e}cembre 2000.}}

\begin{document}
\selectlanguage{french}
\maketitle
\noindent\makebox[2cm][l]{\bf R{\'e}sum{\'e}.}\parbox[t]{11.5cm}{\footnotesize On d{\'e}crit dans cette note une structure
  de torseur sur 
  le sch{\'e}ma affine d{\'e}fini par une collection de relations
  $\QM$-alg{\'e}briques entre les polyz{\^e}tas aux racines de l'unit{\'e}, \ie{}\
  les valeurs prises sur un groupe fix{\'e} de racines complexes de
  l'unit{\'e} par les fonctions hyperlogarithmiques {\`a} plusieurs
  variables. Dans le cas o{\`u} la seule racine de l'unit{\'e}
  est 1 (cas des {\em polyz{\^e}tas}), ces relations sont cens{\'e}es {\^e}tre
  les seules entre ces nombres et le torseur obtenu devrait {\^e}tre {\'e}gal
  {\`a} celui des associateurs de Drinfel'd.  Les formules utilis{\'e}es sont
  en g{\'e}n{\'e}ral celles   de l'action de $\gal(\ov\QM/\QM)$ sur 
  la droite projective priv{\'e}e d'un nombre fini de points. }
\vskip 10pt\selectlanguage{english}
\noindent\hspace{2cm}\parbox{11.5cm}{\bf\em Torsor structure
  associated to some  algebraic relations between
  polyzetas at roots of unity}
\vskip 10pt
\noindent\makebox[2cm][l]{\bf\em Abstract.~}\parbox[t]{11.5cm}{\footnotesize\em
We describe in this note a torsor structure arising on the affine
  scheme defined by a system of $\QM$-algebraic relations between
  polyzetas at roots of unity, \ie{} values of hyperlogarithmic
  functions on a fixed finite subgroup $\G$ of $\CM$. For $\G=\{1\}$,
  ({\em polyzetas} case),
  these relations are believed to span all the $\QM$-algebraic
  relations between those numbers and this torsor should be equal to
  Drinfeld's, \ie{} the Grothendieck-Teichm{\"u}ller group acting on
  associators. The formulas for the general case are derived from the
  action of $\gal(\ov\QM/\QM)$ on the projective line minus finitely
  many points.}

\vspace{3em}

\selectlanguage{english}

\section{Abridged English Version}
In this note, we study some $\QM$-algebraic relations between
some values of the following hyperlogarithmic functions
\begin{equation}\label{def_hyplog}
  \Li_{s_1,\ldots,s_r}(z_1,\ldots,z_r) =
  \sum_{n_1>n_2>\cdots>n_r>0}\frac{z_1^{n_1}z_2^{n_2}\cdots
  z_r^{n_r}}{n_1^{s_1}n_2^{s_2}\cdots n_r^{s_r}},
\end{equation}
where $s_1,\ldots,s_r$ are positive integers. The radius of
convergence of this power series is 1. It is convergent on the
polycircle, except for $(s_1,z_1)=(1,1)$. We shall consider the
numbers $\Li_{s_1,\ldots,s_r}(\s_1,\ldots,\s_r)$, where
$\s_1,\ldots,\s_r$ lie in some  finite multiplicative
subgroup $\Gamma$ of $\CM^*$. We call $r$ the {\em length} of
(\ref{def_hyplog}) and 
$s_1+\cdots+s_r$ its {\em weight}. For $\Gamma=\{1\}$, these are the {\em 
  polyzetas} numbers (also called MZVs, multizetas, \dots) on which we
focused in our thesis (\cite{thesejoe}).  All the
statements present in this note are proved in full detail in
\cite{thesejoe} in this particular case, from which the general one is
not very different. 
  
  In sections \ref{sec_M1} through \ref{sec_D}, we recall
  some fundamental $\QM$-algebraic   relations that those numbers do
  satisfy. For $\G=\{1\}$, it is 
  believed, after a lot of computer checking (\cf{}\ \cite{PetMinh}),
  that these relations span {\em all} the $\QM$-algebraic
  relations between polyzetas. For the general case, the situation
  seems to be more complicated. We will denote by DMRD (Double
  ``M{\'e}lange''\footnote{This is the french word for ``shuffle''.},
  Regularisation and Distribution) the relations given 
  in sections \ref{sec_M1}--\ref{sec_D}.

  As this work grew out of the striking connection between polyzetas
  and associators, let us recall some facts about these. Drinfel'd
  (\cite{DrinQTQH}) introduced the concept of 
associator as a 
step towards the construction of quantized universal enveloping algebras. 
An associator with values in some $\QM$-ring $\k$ is a Lie exponential in
$\sernc{\k}{A,B}$ satisfying some relations (inherited from McLane's
coherence constraints). Among them is the hexagonal
equation, which depends on a parameter $\l\in \k$. Let $\Ass(\k)$ be the
set of all associators and $\Ass_\l(\k)$ the set of
associators with parameter $\l$.  Drinfel'd proved the existence of a
  pro-unipotent group scheme
$\GRT_1(\k)$ acting freely and transitively on each $\Ass_\l(\k)$. This
 is his ``graded''  prounipotent version of the
Grothendieck-Teichm{\"u}ller group,  an interesting object of study
in itself, because of its relations with the absolute Galois group
  (\cf{} \cite{DrinQTQH,IhICM}). 

The only explicitly defined example of an associator, denoted
$\Phi_\KZ$, is obtained in 
terms of solutions of the Knizhnik-Zamolodchikov differential system $KZ_3$
which is also closely related with the iterated integrals that give rise to
the first shuffle relations. It is an element of $\Ass_{i\pi}(\CM)$.  
In 1993 Le and Murakami gave an expression of $\Phi_\KZ$ involving all
the polyzetas. Actually, $\Phi_\KZ$ is almost equal 
to the non-commutative generating series we shall denote
$\LC_\ssh(\triv)$ in section \ref{sec_M1}. As the relations DMRD 
should span all the $\QM$-algebraic relations between polyzetas, by
replacing in $\Phi_\KZ$ the polyzetas by formal symbols satisfying
these relations, we should get a ``universal'' associator over $\QM$,
but we are not able to prove that. 

As those relations are at their simplest on the generating series, we
will consider all the series, with coefficients in an arbitrary
$\QM$-ring $\k$ satisfying them. 
For arbitrary $\G$, this gives a functor
$\k\mapsto\DMRD(\G)(\k)$. We consider it as a concrete realisation of the
affine scheme associated with the relations DMRD. 
The main
result, given in section \ref{sec_tors}, is that this object is
fibered over the affine line $\AM^1$ and is a trivial torsor (as an affine
$\AM^1$-scheme). 
It is indeed very close to Drinfeld's : for $\G=\{1\}$, the torsors
$\DMRD(\G)$ and $\Ass$ work with the same formula.  
The law $\pmt$ we use for the action in 
the general case is not original either : as for $\G=\{1\}$, it comes
from the theory of 
Galois representation into pro-algebraic fundamental groups and from
monodromy torsors (\cf{} \cite{GonchDuke, Wojt}). 
The law $\pmt$ defines a 
pro-unipotent group scheme on its maximal domain of definition. It
is denoted $\MT(\G)$ and described in section \ref{sec_MT}. 
Its infinitesimal structure is
described by two kinds of operators~: the special derivations, which are
relevant for its Lie algebra $\mt(\G)$, and the tangential action
operators (denoted $s_\psi$), which are relevant for the exponential
map $\mt(\G)\to\MT(\G)$.  

The most difficult part of this work is to show the stability of
the second shuffle relation (section \ref{sec_M2}) under this action, because
we're not able to explain it easily in terms of Galois actions or
monodromy actions. The end of this note is devoted to a sketch of 
proof for the main result. 

In section \ref{sec_dm}, we study the tangent space 
$\dmrd(\G)$ near 1 of $\DMRD$ and define a linear form
$\sernc{\k}{X_\G}(\k)\to\k$ who will eventually become the map of the
theorem. We denote by $\dmrd_0(\G)(\k)$ its kernel in
$\dmrd(\G)(\k)$. Next, in section \ref{sec_tangente}, we sketch the proof of
the hardest part of the theorem~: $\DMRD_\l(\G)(\k)$ is stable under the
exponential of the tangential action of $\dmrd_0(\G)(\k)$ on
$\sernc{\k}{X_\G}$ with respect to the law $\pmt$. The arguments used there
prove also that $\dmrd_0(\G)(\k)$ is a Lie subalgebra of $\mt(\G)$,
and therefore its exponential is a sub-group scheme of $\MT(\G)$. 
The section \ref{sec_trans} is devoted to  the transitivity of this
action, thus concluding the proof of theorem. The method we use there
is valid in a much more general framework. 

In the last section, we explain some consequences of the torsor
structure. First, the formal algebra $\DC\MC\RC\DC(\G)$ is a
polynomial algebra and every basis of the finite part of $\dmrd_0(\G)$
gives rise to a set of free generators. For $\G=\{1\}$, this is
{\'E}calle's theorem (\cite{Ecallum}). Next we show that Drinfeld's
irreducible elements, which conjecturally span freely the Lie algebra
of $\GRT_1$, are in $\dmrd_0(\{1\})$, thus enforcing the conjecture that
the relations 
DMRD are actually equivalent to the equations defining associators.  
\subsection{Differences between versions 2 and 3}
  \begin{itemize}
    \item Bad coefficients in projections used for the formulation of
      the distribution relations (section \ref{sec_D}). Consequently,
      all distribution-like formulas have been corrected. This doesn't
      change anything on the result. 
    \item The short exact sequence used in definition 8.1 was
      irrelevant. Also, the coefficients appearing in weight 1 relations change
      with the imbedding of $\G$ into $\CM^*$. So $\DMRD(\G)$ doesn't
      make sense for ``abstract'' $\G$. 
    \item Some tables of dimensions have been included for the tangent spaces 
$\dmrd(\G)$ and $\dmr(\G)$ in which the relevant Lie algebras are of codimension 1. 
   \end{itemize} 

\selectlanguage{french}
\section{Remerciement}
Je suis tr{\`e}s reconnaissant envers Pierre Cartier de m'avoir donn{\'e} ce
sujet de th{\`e}se. Son aide m'a {\'e}t{\'e} pr{\'e}cieuse tant par l'orientation
g{\'e}n{\'e}rale que par l'\oe uvre de simplification et reformulation qu'il a
accomplie, sans laquelle ce travail n'aurait pu {\^e}tre fait. Merci
{\'e}galement {\`a} Michel Petitot pour ses tables de relations et {\`a}
Benjamin Enriquez et Pierre Deligne pour leurs remarques.

\section{Conventions, notations}\label{sec_conv}
Un $\QM$-anneau est un anneau commutatif qui contient
$\QM$. Toutes les alg{\`e}bres sont associatives et unif{\`e}res. Pour tout
$\QM$-anneau $\k$ et tout alphabet $Z$, on note respectivement $Z^*$,
$\assl{\k}{Z}$ et $\sernc{\k}{Z}$ le mono{\"\i}de libre et les $\k$-alg{\`e}bres
de polyn{\^o}mes et de s{\'e}ries non-commutatives form{\'e}es sur $Z$. Dans les
situations consid{\'e}r{\'e}es on aura sur  $\assl{\k}{Z}$ une graduation (le
poids) permetttant les consid{\'e}rations topologiques usuelles. On note
alors $(\cdot|\cdot)$ le produit scalaire de $\assl{\k}{Z}$ pour
lequel $Z^*$ est orthonormale. Cela permettra toujours d'identifier
$\assl{\k}{Z}$ {\`a} son dual gradu{\'e}. Le symbole $\whot$ d{\'e}signe un
produit tensoriel compl{\'e}t{\'e}. On note $1$ le mot vide et $\triv$ le
groupe {\`a} un {\'e}l{\'e}ment. On consid{\`e}re les sch{\'e}mas sur $\spec\QM$ comme des
foncteurs de la cat{\'e}gorie des $\QM$-anneaux dans celle des ensembles
(\cf{} \cite{Demaz}).

\section{Premi{\`e}res relations de m{\'e}lange}\label{sec_M1}
Le contenu de cette partie est tr{\`e}s bien connu. La lettre $\G$ d{\'e}signe
un sous-groupe multiplicatif fini de $\CM^*$. Soit
$X_\G=\{x_0\}\cup\{(x_\s)_{\s\in\G}\}$ un alphabet. On 
munit $\assl{\QM}{X_\G}$ du coproduit $\D$ qui en fait la big{\`e}bre
enveloppante de l'alg{\`e}bre de Lie libre form{\'e}e sur $X_\G$. On la munit
de la graduation pour laquelle les {\'e}l{\'e}ments de $X_\G$ sont de degr{\'e}
1. Dans la suite, cette graduation sera appel{\'e}e le {\em poids}.  
Le produit scalaire $(\cdot|\cdot)$ permet alors de consid{\'e}rer  l'alg{\`e}bre
$(\assl{\QM}{X_\G},\sh)$ duale de la cog{\`e}bre 
$(\assl{\QM}{X_\G}, \D)$ (voir \cite{Reut}). 
{\`A} tout mot $w=x_{\eps_1}x_{\eps_2}\cdots x_{\eps_r}$ en $X_\G$, on associe
l'int{\'e}grale it{\'e}r{\'e}e 
\begin{equation}\label{eq_intiter}
 I(w) = \int_{0\leq t_r \leq \cdots\leq t_1\leq 1}
\bigwedge_{i=1}^r\omega_{\eps_i}(t_i),\qtext{avec}
\omega_0=\frac{dt}{t}\qtext{et} \omega_\s=\frac{\s dt}{1-\s t},\
\text{pour tout}\ \s\in\G
\end{equation}
On pose par convention $I(1)=1$. L'int{\'e}grale $I(w)$ est convergente si
et seulement si le mot 
$w$ ne commence pas par $x_1$ et ne se finit pas par $x_0$. On note
$X_{\G,\cv}^*$ l'ensemble de ces mots et $\QM(X_{\G,\cv}^*)$ le
sous-$\QM$-espace vectoriel de $\assl{\QM}{X_\G}$ qu'il engendre. 
On {\'e}tend l'application $I$ par lin{\'e}arit{\'e} {\`a} $\QM(X_{\G,\cv}^*)$. 
\begin{prop}
   L'ensemble $\QM(X_{\G,\cv}^*)$ est une sous-alg{\`e}bre de
   $(\assl{\QM}{X_\G},\sh)$. L'application $I$ est un morphisme
   d'alg{\`e}bres de  $(\QM(X_{\G,\cv}^*),\sh)$ dans $\CM$. 
\end{prop} 
Ceci constitue un premier syst{\`e}me de relations entre polyz{\^e}tas aux
racines de l'unit{\'e}, gr{\^a}ce {\`a} la proposition ci-dessous
\begin{prop}
   Pour tous entiers strictement positifs $s_1,\ldots,s_r$ et tous
   {\'e}l{\'e}ments $\s_1,\ldots,\s_r$ de $\G$ tels que $(s_1,\s_1)$ soit
   diff{\'e}rent de $(1,1)$, on a~:
\begin{equation}\label{eq_corr}
 \Li_{s_1,\ldots,s_r}(\s_1,\ldots,\s_r) =
I(x_0^{s_1-1}x_{\s_1}x_0^{s_2-1}x_{\s_1\s_2}\cdots
x_0^{s_r-1}x_{\s_1\cdots\s_r}) 
\end{equation}
\end{prop}
Il est plus satisfaisant d'{\'e}tendre la d{\'e}finition de $I$ {\`a}
$\assl{\QM}{X_\G}$ tout entier. Cette {\em r{\'e}gularisation} est rendue
possible gr{\^a}ce aux propri{\'e}t{\'e}s du produit $\sh$~:
\begin{prop}\label{prop_regul}
  Il existe un unique morphisme d'alg{\`e}bres $\ov{I}$ de $(\assl{\QM}{X_\G},\sh)$
  dans $\CM$ co{\"\i}ncidant avec $I$ sur $\QM(X_{\G,\cv}^*)$ et v{\'e}rifiant
  $\ov{I}(x_0)=\ov{I}(x_1)=0$.
\end{prop} 
De fa\c con duale, consid{\'e}rons la s{\'e}rie g{\'e}n{\'e}ratrice
non-commutative $\LC_\ssh(\G)\ \ass\ \sum_{w\in X_\G^*} \ov{I}(w)w$ qui
appartient {\`a} $\sernc{\CM}{X_\G}$. Par
dualit{\'e} entre le produit $\sh$ et le coproduit $\D$, la premi{\`e}re
relation de m{\'e}lange est {\'e}quivalente {\`a} 
$\D\LC_\ssh(\G)=\LC_\ssh(\G)\whot\LC_\ssh(\G)$. Si l'on identifie
$A$ avec $x_0$ et $B$ avec $-x_1$, l'associateur $\Phi_\KZ$ est {\'e}gal {\`a}
$\LC_\ssh(\triv)$. 
\section{Secondes relations de m{\'e}lange}\label{sec_M2}
Il s'agit ici de d{\'e}crire la combinatoire des relations que l'on
obtient par produit de s{\'e}ries et d{\'e}composition du domaine de
sommation, comme par exemple~:
\begin{eqnarray}
\Li_{s_1}(\s_1)\Li_{s_2}(\s_2)&=&
\sum_{n_1,n_2>0}\frac{\s_1^{n_1}\s_2^{n_2}}{n_1^{s_1}n_2^{s_2}}=
\sum_{n_1>n_2>0}\frac{\s_1^{n_1}\s_2^{n_2}}{n_1^{s_1}n_2^{s_2}}+
\sum_{n_2>n_1>0}\frac{\s_1^{n_1}\s_2^{n_2}}{n_1^{s_1}n_2^{s_2}}+
\sum_{n>0}\frac{(\s_1\s_2)^n}{n^{s_1+s_2}}\nonumber\\&=&\Li_{s_1,s_2}(\s_1,\s_2)
+ \Li_{s_2,s_1}(\s_2,\s_1)+\Li_{s_1+s_2}(\s_1\s_2)\label{exemple}
\end{eqnarray}
Dans le cas $\G=\triv$, on sait que l'on peut d{\'e}crire ces relations
en interpr{\'e}tant $\Li$ comme un op{\'e}rateur de sp{\'e}cialisation sur l'alg{\`e}bre
des fonctions quasi-sym{\'e}triques. Dans \cite{thesejoe}, on se ram{\`e}ne {\`a}
la cog{\`e}bre duale. Dans \cite{Big}, il est propos{\'e} une g{\'e}n{\'e}ralisation~:
les \og fonctions quasi-sym{\'e}triques color{\'e}es\fg{}. On peut {\'e}galement
la d{\'e}crire en termes de cog{\`e}bre~: 

Soit $Y_\G=\{y_{n,\nu}\}_{n\in\NM^*, \nu\in\Gamma}$ un alphabet. On
munit $\assl{\QM}{Y_\G}$ de la graduation telle que $y_{n,\nu}$ soit
de degr{\'e} $n$ pour tous $n$ et $\nu$. On l'appelle encore le poids. 
L'unique morphisme d'alg{\`e}bres de $\assl{\QM}{Y_\G}$ dans son carr{\'e}
tensoriel v{\'e}rifiant~:
\begin{equation}
\Det y_{n,\nu} = \underset{k,l\in\NM,\ \kappa,\l\in\G}{\sum_{k+l=n,
    \kappa\l=\nu}} y_{k,\kappa}\ot y_{l,\lambda}\qtext{avec la
    convention} y_{0,\s}=\left\{\begin{array}{lcl} 
1 &\text{si}& \s=1\\
0 &\text{si}& \s\neq1
\end{array}\right.
\end{equation}
fait de $(\assl{\QM}{Y_\G},\cdot,\Det)$ une big{\`e}bre gradu{\'e}e. On note $\et$ le
produit dual du coproduit $\Det$. Soit $Y_{\G,\cv}^*$ l'ensemble des
mots de $Y_\G^\et$ ne commen\c cant pas par $y_{1,1}$ et
$\QM(Y_{\G,\cv}^*)$ le sous-$\QM$-espace vectoriel qu'il engendre.  On
note, par abus, encore $\Li$ l'application $\QM$-lin{\'e}aire de
$\QM(Y_{\G,\cv}^*)$ dans $\CM$ qui v{\'e}rifie~:
$$ L(y_{s_1,\s_1}y_{s_2,\s_2}\cdots y_{s_r,\s_r}) =
\Li_{s_1,s_2,\ldots,s_r}(\s_1,\s_2,\ldots,\s_r)\qtext{et}L(1)=1$$
\begin{prop}[Seconde relation de m{\'e}lange]
   L'ensemble $\QM(Y_{\G,\cv}^*)$ est une sous-alg{\`e}bre de
   $(\assl{\QM}{Y_\G},\et)$ et l'application lin{\'e}aire $\Li$ est un morphisme
   d'alg{\`e}bres de $(\QM(Y_{\G,\cv}^*),\et)$ dans $\CM$. 
\end{prop}
Ceci exprime les relations du type (\ref{exemple}). Dans le cas
$\G=\triv$, il est classique que la big{\`e}bre
$(\assl{\QM}{Y_\G},\cdot,\Det)$ se ram{\`e}ne {\`a} une big{\`e}bre enveloppante
d'alg{\`e}bre de Lie libre (\cf{} \cite{MalReut}). La g{\'e}n{\'e}ralisation au
cas $\G$ quelconque est 
facile et permet d'{\'e}tendre $\Li$ {\`a} $\assl{\QM}{Y_\G}$, comme {\`a} la
section pr{\'e}c{\'e}dente~: 
\begin{prop}
   Il existe une unique s{\'e}rie $\LC_\et(\G)\in\sernc{\CM}{Y_\G}$ telle
   que 
$$(\LC_\et(\G)|y_{1,1})=0,\quad (\LC_\et(\G)|w)=L(w)\ \forall w\in
Y_{\G,\cv}^*\qtext{et} \Det(\LC_\et(\G))=\LC_\et(\G)\whot\LC_\et(\G)$$
\end{prop}

\section{Relations de r{\'e}gularisation}\label{sec_R}
Le codage dans l'alphabet $Y$ que l'on vient d'utiliser est quasiment
tautologique. Pour tenir compte de (\ref{eq_corr}) il est commode 
d'identifier, pour tout $(n,\nu)$ appartenant {\`a} $\NM^*\times\G$,
la lettre $y_{n,\nu}$ au mot $x_0^{n-1}x_\nu$. Cela permet de plonger
par morphisme 
pour le produit de concat{\'e}nation $\assl{\QM}{Y_\G}$ dans
$\assl{\QM}{X_\G}$. On notera $\pi_Y :
\sernc{\QM}{X_\G}\to\sernc{\QM}{Y_\G}$ la projection duale de ce
plongement.
On constate que $\QM(Y_{\G,\cv}^*)$ est {\'e}gal {\`a}
$\QM(X_{\G,\cv}^*)$. Si l'on note $\ps$ (produits successifs)
l'endomorphisme lin{\'e}aire de $\assl{\QM}{Y_\G}$ donn{\'e} par 
\begin{equation}
\forall s_1,\ldots,s_r \in\NM^*,\ \s_1,\ldots,\s_r\in\G,\
\ps(y_{s_1,\s_1}y_{s_2,\s_2}\cdots y_{s_r,\s_r}) =
y_{s_1,\s_1}y_{s_2,\s_1\s_2} \cdots y_{s_r,\s_1\s_2\cdots\s_r}
\end{equation}
la formule (\ref{eq_corr}) devient alors~:  $\forall v\in\QM(Y^*_{\G,\cv}),
L(w)=I(\ps(w))$. La question se pose de savoir si cette formule reste
vraie pour $\Li_\et$ et $I_\ssh$. Ce serait {\'e}quivalent, en notant
$\qs$ (quotients successifs) l'inverse de $\ps$, {\`a}
$\qs\pi_Y(\LC_\ssh(\G))=\LC_\et(\G)$. Ceci est malheureusement faux,
mais {\'E}calle (\cite{Ecalle}) a donn{\'e}  l'expression g{\'e}n{\'e}rale de la
correction {\`a} apporter~: 
\begin{prop}[Relation de r{\'e}gularisation] On a 
$$ \LC_\et = \exp\left(\sum_{n\geq 2} \frac{(-1)^{n-1}}{n}
  \z(n)y_1^n\right) \qs\pi_Y(\LC_\sh)$$
\end{prop}
Cette formule implique la relation d'Hoffman utilis{\'e}e en
\cite{PetMinh} et g{\'e}n{\'e}ralis{\'e}e {\`a} $\G$ quelconque en \cite{Big}~:
celle-ci se r{\'e}duit {\`a} l'absence de terme $n=1$ dans la somme ci-dessus.
Faute de disposer de  la d{\'e}monstration d'{\'E}calle, on en a donn{\'e} une autre
dans \cite{thesejoe} pour le cas $\G=\triv$, bas{\'e}e sur les remarques de
Boutet de Monvel (\cite{Boutet}). La g{\'e}n{\'e}ralisation {\`a} $\G$ quelconque
en est imm{\'e}diate. 

\section{Relations de distribution}\label{sec_D}
On utilise une description proche d'une de celles de Goncharov
(\cf{} \cite{GonchDuke}). Pour tout sous-groupe $\G'$ de $\G$, soit
$\proj^1_{\G\to\G'}$ la projection de $\sernc{\CM}{X_\G}$ sur
$\sernc{\CM}{X_{\G'}}$ duale de l'inclusion de $\assl{\CM}{X_{\G'}}$
dans $\assl{\CM}{X_\G}$ et $\proj^2_{\Gamma\to\Gamma'}$ le
morphisme d'alg{\`e}bres topologiques de $\sernc{\CM}{X_\G}$ dans
$\sernc{\CM}{X_{\G'}}$ caract{\'e}ris{\'e} par 
\begin{eqnarray*}
\forall\s\in\G,\quad  \proj^2_{\Gamma\to\Gamma'}(x_\s) &=& x_{\s^{[\G:\G']}}\\
\qtext{et}\proj^2_{\Gamma\to\Gamma'}(x_0) &=& [\G:\G']x_0
\end{eqnarray*} 
Avec ces notations, on a~:
\begin{equation}
 \proj^2_{\Gamma\to\Gamma'}(\LC_\ssh(\G)) =
\exp\left(\sum_{\s^{[\G:\G']}=1,\s\neq1}\Li_1(\s)x_1\right)\proj^1_{\Gamma\to\Gamma'}(\LC_\ssh(\G)), 
\end{equation}

Pour {\^e}tre plus pr{\'e}cis, les relations de distribution expriment l'{\'e}galit{\'e}
des coefficients des mots convergents dans $\proj^1_{\G\to\G'}(\LC_\ssh(\G))$ et
$\proj^2_{\G\to\G'}(\LC_\ssh(\G))$. Ces deux s{\'e}ries sont
\og{}group-like\fg{} dans $(\sernc{\CM}{X_{\G'}},\D)$, car les deux
projections sont des morphismes de cog{\`e}bres. D'apr{\`e}s une variante plus
g{\'e}n{\'e}rale de la proposition \ref{prop_regul}, il suffit donc
pour comparer ces s{\'e}ries d'{\'e}tudier leurs termes en $x_1$, ce qui donne
la formule ci-dessus.  

Enfin, il faut ajouter {\`a} toutes ces relations celles que l'on obtient
en exprimant que, pour tout $\s\in\G$, les
$\log(1-\s)-\log(1-\s^{-1})$ sont des multiples 
rationnels de $i\pi$ et donc proportionnels sur $\QM$. On les
appellera {\em relations de poids 1} et on les inclura dans le syst{\`e}me
DMRD. 
\section{Objets formels}\label{sec_defDM}
Pour {\'e}tudier la structure alg{\'e}brique sous-jacente aux relations
DMRD, on pourrait consid{\'e}rer la 
$\QM$-alg{\`e}bre $\DC\MC\RC\DC(\Gamma)$
engendr{\'e}e par des symboles formels 
index{\'e}s comme les 
valeurs des hyperlogarithmes aux racines de l'unit{\'e}, modulo les
relations DMRD. Il est {\'e}quivalent, et plus simple, de travailler sur
l'objet dual, \ie{} l'ensemble des 
s{\'e}ries satisfaisant aux m{\^e}mes relations que $\LC_\ssh(\G)$, mais {\`a}
coefficients dans des $\QM$-anneaux quelconques. 

\begin{defi}\label{def_DM}
   Pour tout $\QM$-anneau $\k$ et tout groupe commutatif fini
   $\Gamma$, on note $\DMR(\Gamma)(\k)$ l'ensemble des {\'e}l{\'e}ments $\Phi$
   de $\sernc{\k}{X_{\G}}$ tels que~:
\begin{eqnarray}
  (\Phi|1) =  1 &\text{\rm et}& (\Phi|x_0) = (\Phi|x_1) = 0
  \label{eq_Phipoids1}\\ 
  \D\Phi\ =\ \Phi\whot_{\k}\Phi &\text{\rm et}&  \Det\Phi_\et\ =\
  \Phi_\et\whot_{\k}\Phi_\et, \\
\qtext{o{\`u} l'on pose} \Phi_\et\ \ass\ \Phi_\corr\cdot\qs\pi_Y(\Phi)&\text{\rm et}&\Phi_\corr\ \ass
 \   \exp\left(\sum_{n\geq2}
{\frac{(-1)^{n-1}}{n}(\pi_Y(\Phi)|y_n)y_1^n}\right)
\end{eqnarray}
Si, de plus, $\G$ est un sous-groupe multiplicatif de $\CM^*$, on note
 $\DMRD(\G)(\k)$ l'ensemble des {\'e}l{\'e}ments $\Phi$ de 
 $\DMR(\G)(\k)$ v{\'e}rifiant la relation de distribution
$$  \proj^2_{\Gamma\to\Gamma'}(\Phi) =
\exp\left(\sum_{\s^{[\G:\G']}=1}(\Phi|x_\s)x_1\right)\proj^1_{\Gamma\to\Gamma'}(\Phi), $$
pour tout sous-groupe $\G'$ de $\G$ et les relations de
poids 1.
\end{defi}
Le symbole $\DMRD(\G)$ est un foncteur
de la cat{\'e}gorie des $\QM$-anneaux dans celle des
ensembles. Il est naturellement isomorphe {\`a}
$\spec{\DC\MC\RC\DC(\G)}$. La s{\'e}rie $\LC_\ssh(\G)$ appartient {\`a}
$\DMRD(\G)(\CM)$ et on a $\LC_\ssh(\triv)=\Phi_\KZ$. Le th{\'e}or{\`e}me
d'{\'E}calle (\cite{Ecalle}), {\'e}nonc{\'e} fin 1999, dit que 
$\DC\MC\RC\DC(\triv)$ est une alg{\`e}bre de polyn{\^o}mes et en d{\'e}crit les
g{\'e}n{\'e}rateurs libres. 

\section{Le groupe $\MT(\G)$ et sa structure infinit{\'e}simale}\label{sec_MT}
Dans ce qui suit, pour tout {\'e}l{\'e}ment $\g$ de $\G$, on note $t_\g$
l'action naturelle de $\G$ sur $\sernc{\k}{X_\G}$, \ie{} 
l'automorphisme de $\k$-alg{\`e}bres topologiques de $\sernc{\k}{X_\G}$ fixant
$x_0$ et envoyant $x_\s$ sur $x_{\g\s}$, pour tout $\s\in\G$. 

Pour tout $\QM$-anneau $\k$, soit $\MT(X_\G)(\k)$ l'ensemble
des s{\'e}ries de
$\sernc{\k}{X_\G}$ dont 
le terme constant vaut 1 et muni du produit $\pmt$ donn{\'e} par $G\pmt H
= G\cdot\kappa_G(H)$,  
o{\`u} $\kappa_G$ est le morphisme de $\k$-alg{\`e}bres topologiques
caract{\'e}ris{\'e} par 
\begin{equation}
   \kappa_G(x_0) = x_0 \qtext{et} \kappa_G(x_\s) = t_\s(G)^{-1}x_\s
   t_\s(G) (\forall\s\in\G)
\end{equation}
Vu comme foncteur en $\k$, c'est un sch{\'e}ma en groupes
pro-unipotent sur $\QM$, car $\kappa$ en est une repr{\'e}sentation
lin{\'e}aire fid{\`e}le dans $\sernc{\k}{X_\G}$ et $\kappa_G(w)=w$ modulo des termes
de plus haut poids que $w$, pour tout mot $w$ de $X_\G^*$. 
Ce groupe n'est autre que celui des automorphismes de
$\sernc{\k}{X_\G}$, {\'e}quivariants pour l'action naturelle de $\G$, qui laissent
stables {\`a} conjugaison pr{\`e}s les $(x_\s)_{\s\in\G}$ et laissent $x_0$
invariant. Son alg{\`e}bre de Lie $\k\mapsto\mt(\G)(\k)$
est form{\'e}e des s{\'e}ries sans terme constant de $\sernc{\k}{X_\G}$. On
note $\psi_1,\psi_2\mapsto\ih{\psi_1}{\psi_2}$ son crochet. En
lin{\'e}arisant le produit $\pmt$ au voisinage de 1, on obtient les
r{\'e}sultats suivants~: 

Pour toute s{\'e}rie $\psi\in\mt(\G)(\k)$, soit $d_\psi$ la d{\'e}rivation
continue de $\sernc{\k}{X_\G}$ caract{\'e}ris{\'e}e par 
\begin{equation}
d_\psi(x_0) = 0\qtext{et} \forall \s\in\G,\ d_\psi(x_\s) = [x_\s, t_\s(\psi)]
\end{equation}
C'est ce que Goncharov, suivant Ihara, appelle une d{\'e}rivation sp{\'e}ciale
{\'e}quivariante (pour l'action de $\G$). 
Si l'on note $s_\psi(\phi)=\psi\phi+d_\psi(\phi)$, pour tout $\phi$ de
$\sernc{\k}{X_\G}$, et $\Exp$ l'application exponentielle de
$\mt(\G)(\k)$ dans $\MT(\G)(\k)$, on a~:
\begin{eqnarray}
\forall\psi\in\mt(\G)(\k),\ H\in\sernc{\k}{X_\G}\ \Exp(\psi)\pmt H &=&
\exp(s_\psi)(H) \\
\forall\psi_1,\psi_2\in\mt(\G)(\k),\ s_{\ih{\psi_1}{\psi_2}} &=&
  [s_{\psi_1},s_{\psi_2}] 
\end{eqnarray} 
De plus, pour tout $\psi$ de $\sernc{\k}{X_\G}$, on a $s_\psi(1)=\psi$. 

\section{{\'E}nonc{\'e} du r{\'e}sultat principal}\label{sec_tors} 

\begin{theo}\label{theo_action}
Si $\Gamma$ est un sous-groupe fini de $\CM$, il existe un
morphisme de sch{\'e}mas $\DMRD(\G)\to\AM^1$ ayant les propri{\'e}t{\'e}s suivantes~:
\begin{itemize}
  \item pour tout $\QM$-anneau $\k$, l'application $\DMRD(\G)(\k)\to\k$
    est surjective
  \item la fibre sp{\'e}ciale $\DMRD_0(\G)$ au-dessus de 0 est un sous-sch{\'e}ma
    en groupes de $\MT(\Gamma)$
  \item pour tout $\QM$-anneau $\k$ et tout $\l\in\k$, le groupe
    $\DMRD_0(\G)(\k)$ agit librement et transitivement sur
    $\DMRD_\l(\k)$ par translation {\`a} gauche (au sens de $\pmt$)
\end{itemize}
\end{theo}

Pour $\G=\triv$ et $\G=\{-1, 1\}$ la fl{\`e}che du th{\'e}or{\`e}me donne la
valeur $\l$ de $\k$ par laquelle on a remplac{\'e} $\z(2)$. Dans le cas o{\`u}
le cardinal de $\G$ vaut au moins 3, elle donne la valeur par laquelle on
a remplac{\'e} $\log(1 - e^{2i\pi/N})-\log(1 - e^{-2i\pi/N})$, lui-m{\^e}me
    multiple rationnel de $i\pi$. 

{\`A} titre de comparaison, on peut r{\'e}sumer ainsi les propositions 5.5 et
5.9 de \cite{DrinQTQH} dont le th{\'e}or{\`e}me \ref{theo_action} est inspir{\'e}~:
\begin{theo}[Drinfel'd]
   Le sch{\'e}ma $\Ass_0$ est un sous-sch{\'e}ma en groupes de $\MT(\triv)$ et
   $\Ass_0(\k)$, {\'e}gal {\`a} $\GRT_1(\k)$, agit librement et transitivement
   par translation {\`a} gauche au sein de $\MT(\triv)(\k)$ sur chaque $\Ass_\l(\k)$, pour tout
   $\QM$-anneau $\k$ et tout $\l\in\k$.  
\end{theo}
La fl{\`e}che $\Ass\to\AM^1$ indique quel est l'{\'e}l{\'e}ment $\l$ de $\k$ qui
intervient dans l'{\'e}quation hexagonale. Dans le cas de $\Phi_\KZ$,
on a $\l=i\pi$.
\begin{conj}\label{conj}
   Pour tout $\QM$-anneau $\k$ et tout $\l\in\k$, on a
   $\Ass_\l(\k)=\DMRD_{-\l^2/6}(\triv)(\k)$.
\end{conj}
Cette conjecture est {\'e}tay{\'e}e par la compatibilit{\'e} entre les conjectures
de dimension de Zagier (\cf{} \cite{ZagierECM}) et la variante de
Drinfel'd de la conjecture de Deligne (\cf{} \cite{DrinQTQH}, 
  p. 860). On donnera quelques informations suppl{\'e}mentaires {\`a} ce
  propos {\`a} la fin de la section \ref{sec_conseq}.

\section{Les espaces tangents $\dmrd$ et $\dmrd_0$}\label{sec_dm}
On consid{\`e}re ici la lin{\'e}arisation au voisinage de 1 de $\DMRD(\G)$~:
\begin{defi}
Pour tout $\QM$-anneau $\k$, soit $\dmrd(\k)$ l'ensemble des s{\'e}ries
$\psi$ de $\sernc{\k}{X_\G}$ qui v{\'e}rifient~:
\begin{eqnarray}
(\psi|x_0)&=&(\psi|x_1)=0 \\
\D\psi=1\ot_\k\psi + \psi\ot_\k 1 &\text{et}& \Det(\psi_\et)=1\ot_\k\psi_\et +
\psi_\et\ot_\k 1,\label{eq_dmprim}\\
\qtext{o{\`u} l'on pose} \psi_\et\ \ass\ \qs\pi_Y(\psi) + \psi_\corr &\text{et}&
\psi_\corr\ \ass\ \sum_{n\geq 2}
\frac{(-1)^{n-1}}{n} (\psi|y_n)y_1^n, 
\end{eqnarray}
ainsi que la relation de distribution tangente 
$$ 
\proj^2_{\Gamma\to\Gamma'}(\Phi) = \proj^1_{\Gamma\to\Gamma'}(\Phi) +
\sum_{\s^{[\G:\G']}=1}(\Phi|x_\s)x_1,
$$ 
pour tout sous-groupe $\Gamma'$ de $\Gamma$. 
et les relations de poids 1, inchang{\'e}es car lin{\'e}aires. 
\end{defi}

\begin{prop}
   Pour tout $\psi$ de $\dmrd(\k)$, homog{\`e}ne de poids $n\geq3$ et tout
   $\nu\in\G$, on a  
\begin{equation}\label{eq_excep}
   (\psi_\et|y_{n,\nu}) + (-1)^n (\psi_\et|y_{n,\nu^{-1}}) = 0 
\end{equation}
\end{prop} 
 On d{\'e}montre cela par une analyse explicite des contraintes portant sur
les termes de longueur 1 et 2 de $\psi$ impos{\'e}es par les {\'e}quations
\ref{eq_dmprim}. 

On qualifiera d'exceptionnel un {\'e}l{\'e}ment de $\dmrd(\G)(\k)$ ne v{\'e}rifiant
pas l'{\'e}quation (\ref{eq_excep}). Soit $\dmrd_0(\k)$ l'ensemble des
{\'e}l{\'e}ments de $\dmrd(\k)$ qui la v{\'e}rifient. Il y a essentiellement un
seul {\'e}l{\'e}ment exceptionnel pour chaque $\G$~: 
\begin{prop}
  Pour tout $\G$, il existe un {\'e}l{\'e}ment homog{\`e}ne $\alpha_\G$ de
  $\dmrd(\G)(\QM)$ tel 
  que $\dmrd_0(\G)(\k)$ soit le noyau dans $\dmrd(\G)(\k)$ de la forme
  lin{\'e}aire $\psi\mapsto(\psi|\alpha_G)$. 
\end{prop}
Si le cardinal de $\G$ est au moins 3 on peut prendre pour $\alpha_G$
un des $y_{1,\nu}-y_{1,\nu^{-1}}$ avec $\nu^2\neq1$ (gr{\^a}ce aux
  relations de poids 1, toutes les formes lin{\'e}aires ainsi obtenues sont
  colin{\'e}aires sur $\QM$). Si $\G=\triv$ ou $\G=\{\pm 1\}$, on peut
  prendre $\alpha_\G=y_2$. On consid{\`e}re dans la suite $\alpha_\G$
  comme fix{\'e}.

\begin{defi}
  Pour tout $\QM$-anneau $\k$ et tout $\l\in\k$, on note
  $\DMRD_\l(\G)(\k)$ l'ensemble des {\'e}l{\'e}ments $\Phi$ de $\DMRD(\G)(\k)$
  v{\'e}rifiant~: $(\Phi|\alpha_\G)=\l$. 
\end{defi}

On a ainsi d{\'e}crit le morphisme de sch{\'e}mas $\DMRD(\G)\to\AM^1$ du
th{\'e}or{\`e}me. 

\section{Action tangente}\label{sec_tangente}
On donne dans cette partie quelques indications sur la d{\'e}monstration
de l'{\'e}nonc{\'e} ci-dessous qui constitue une moiti{\'e} du th{\'e}or{\`e}me. Pour
all{\'e}ger, les termes \og $\et$-primitif\fg{} et \og
$\et$-cod{\'e}rivation\fg{} signifient respectivement \og primitif\fg{} et
\og cod{\'e}rivation\fg{} pour le coproduit $\Det$. 
\begin{prop}\label{prop_tang}
   Pour tout $\QM$-anneau $\k$, tout $\l\in\k$ et tout
   {\'e}l{\'e}ment $\psi$ de $\dmrd_0(\G)(\k)$, l'ensemble $\DMRD_\l(\G)(\k)$
   est stable 
   par multiplication $\pmt$ {\`a} gauche par $\Exp(\psi)$.
\end{prop}

Soit $\psi\in\dmrd_0(\k)$. La stabilit{\'e} des
{\'e}quations (\ref{eq_Phipoids1}) par $\exp(s_\psi)$ est
{\'e}vidente. 

En ce qui concerne les relations de m{\'e}lange, qui s'expriment pour
$\Phi\in\DMRD(\G)(\k)$ par le fait que $\Phi$ et $\Phi_\et$ sont 
\og{}group-like\fg{} respectivement pour $\D$ et $\Det$, la situation est
dissym{\'e}trique. En effet, 
l'ensemble des {\'e}l{\'e}ments \og{}group-like\fg{} de $(\sernc{\k}{X_\G},\D)$
est stable par $\pmt$. Dans l'optique qui est la notre, cela peut se
formuler ainsi~: soit $\psi\in\sernc{\k}{X_\G}$ ;
l'op{\'e}rateur $s_\psi$ est somme de la d{\'e}rivation $d_\psi$ et de la
multiplication {\`a} gauche par $\psi$ ; si $\psi$ est primitif (pour
$\D$), cette derni{\`e}re est une cod{\'e}rivation, ainsi que $d_\psi$, car c'est
une d{\'e}rivation qui envoie les lettres de $X_\G$ sur des {\'e}l{\'e}ments
primitifs; $\exp(s_\psi)$ est donc un morphisme de cog{\`e}bres, ce qui
permet de conclure. On va transposer cette m{\'e}thode au coproduit
$\Det$, mais on ne peut rien esp{\'e}rer d'aussi simple, car l'ensemble des
\og{}group-like\fg{} pour $\Det$ n'est pas stable par $\exp s_\psi$, si
$\psi$ est un $\et$-primitif quelconque. 

Premi{\`e}rement, pour tout $\psi\in\sernc{\k}{X_\G}$, on voit facilement
que le noyau de $\pi_Y$ est stable par $s_\psi$. On peut donc
consid{\'e}rer l'endomorphisme $\k$-lin{\'e}aire $s_\psi^Y$ de
$\sernc{\k}{Y_\G}$, quotient par $\qs\pi_Y$. D'autre part, si l'on
note $\d_{x_0}$ la 
d{\'e}riv{\'e}e partielle par rapport {\`a} $x_0$ de $\sernc{\k}{X_\G}$, la
projection $\pi_Y$ est bijective de $\ker\d_{x_0}$ sur
$\sernc{\k}{Y_\G}$. Son inverse $\s$ est donn{\'e} par 
\begin{equation}\label{eq_sec}
 \s\psi = \sum_{i\geq0} \frac{(-1)^i}{i!} \d_{x_0}^i(\psi)x_0^i
\end{equation}
Si $\psi$ est une s{\'e}rie de Lie, $\d_{x_0}(\psi)$ et $(\psi|x_0).1$
co{\"\i}ncident. Si $\psi$ appartient {\`a} $\dmrd(\G)(\k)$, on a donc 
$\psi=\s\ps(\qs\pi_Y\psi)$. Comme $s_\psi$ d{\'e}pend lin{\'e}airement de
$\psi$, pour $\psi\in\dmrd(\G)(\k)$, on peut d{\'e}composer $s_\psi$ en
$s_{\s\ps\psi_\et} - s_{\s\psi_\corr}$.
\begin{prop}
   Pour un {\'e}l{\'e}ment $\psi$ de $\dmrd_0(\G)(\k)$, homog{\`e}ne de poids $p$,
   l'endomorphisme $\k$-lin{\'e}aire 
   $s_{\s\ps\psi_\et}$ de $\sernc{\k}{Y_\G}$ est une
   $\et$-cod{\'e}rivation 
\end{prop}
Pour d{\'e}montrer cette proposition, on d{\'e}compose $s_{\s\ps\psi_\et}$ en
somme de la translation {\`a} droite par $\psi_\et$ qui est une
$\et$-cod{\'e}rivation, $\psi_\et$ {\'e}tant $\et$-primitif et d'une
d{\'e}rivation. On exprime les valeurs de celle-ci sur les $Y_{n,\nu}$ en
fonction des $Y_{n,\nu}$, de $\psi_\et$ et d'op{\'e}rateurs respectant la
$\et$-primitivit{\'e}. On teste alors l'identit{\'e} de cod{\'e}rivation sur
les $Y_\G$; elle finit par se ramener {\`a} (\ref{eq_excep}).

L'op{\'e}rateur $\exp s_{\s\psi_\et}^Y$ est donc un automorphisme de 
$\k$-cog{\`e}bres topologiques de $(\sernc{\k}{Y_\G},\Det)$. 
Si $\Phi$ appartient {\`a} $\DMRD_\l(\k)$, l'{\'e}l{\'e}ment $\exp
s_{\s\ps\psi_\et}^Y(\Phi_\et)$ est donc \og{}group-like \fg{} pour $\Det$. 
Pour conclure, il suffit de prouver que c'est exactement
$(\exp(s_\psi)(\Phi))_\et$. En effet, les termes
correctifs (du type $\Phi_\corr$ ou $\psi_\corr$) sont des
s{\'e}ries en $y_1$ et $x_1=y_1$ est central ; de plus, pour tout
$G\in\sernc{\k}{X_\G}$, on a $y_1\pmt G = y_1G$ et
$(\exp(s_\psi)(G)|y_n)=(\psi|y_n)+(G|y_n)$. 
  Pour les relations de distribution, il suffit de voir que les
  coefficients des $(x_\s)_{\s\in\G}$ s'additionnent de la m{\^e}me fa\c
  con et que les applications $\proj^1_{\G\to\G'}$ et
  $\proj^2_{\G\to\G'}$ sont des  morphismes de $\MT(\G)$ dans $\MT(\G')$. 

  La proposition \ref{prop_tang} implique {\'e}galement que
  $\dmrd_0(\G)(\k)$ est une sous-alg{\`e}bre de Lie de $\mt(\G)$, car 
  $\ih{\psi_1}{\psi_2}=[s_{\psi_1},s_{\psi_2}](1)$, pour tous
  $\psi_1,\psi_2$ de $\dmrd(\G)$, ce qui prouve que
  $\ih{\psi_1}{\psi_2}$ v{\'e}rifie les {\'e}quations (\ref{eq_dmprim}). Les
  autres {\'e}quations se traitent comme ci-dessus. Par la formule de Campbell-Hausdorff, $\Exp(\dmrd_0(\k))$ est donc bien un sous-groupe de
  $\MT(\G)$.

\section{Transitivit{\'e}}\label{sec_trans}
On indique ici comment se prouve la deuxi{\`e}me partie du th{\'e}or{\`e}me~:
\begin{prop}
   Pour tout $\QM$-anneau $\k$ et tout $\l\in\k$, l'action de
   $\Exp(\dmrd_0(\k))$ sur $\DMRD_\l(\k)$ est libre et transitive. 
\end{prop}
La libert{\'e} est {\'e}vidente, puisque l'action est une translation au sein
d'un groupe. Pour obtenir la transitivit{\'e}, on suit la m{\'e}thode
d'approximations successives expos{\'e}e dans \cite{BarNatGT}.

Pour tout entier $n$, on note $\assln{\k}{X_\G}{n}$ le quotient de
$\sernc{\k}{X_\G}$ par le $n+1$\eme terme de la filtration associ{\'e}e au
poids et $\tr{\pi}{n}$ la projection correspondante. 
 On consid{\`e}re $\assln{\k}{X_\G}{n}$ comme inclus dans
$\sernc{\k}{X_\G}$ : la projection $\tr{\pi}{n}$ envoie tout mot de
poids au moins $n+1$ sur 0 et laisse les autres fixes. On note
$\tr{\DMRD_\l}{n}(\G)(\k)$ l'ensemble des {\'e}l{\'e}ments de $\assln{\k}{X_\G}{n}$
satisfaisant aux {\'e}quations d{\'e}finissant $\DMRD_\l(\G)(\k)$, modulo des
termes de poids au moins $n+1$. Les $\tr{\DMRD_\l}{n}(\G)(\k)$ forment
un syst{\`e}me projectif dont la limite est $\DMRD_\l(\G)(\k)$. 
 
La remarque principale est la suivante~: {\'e}tant donn{\'e} un {\'e}l{\'e}ment 
$\Phi=\Phi_1+\cdots+\Phi_n$ de $\tr{\DMRD}{n}(\G)(\k)$, si l'on
cherche un {\'e}l{\'e}ment $\Phi_{n+1}$ de $\assl{\k}{X_\G}$, homog{\`e}ne de
poids $n+1$ tel que $\Phi+\Phi_{n+1}$ appartienne {\`a}
$\tr{\DMRD}{n+1}(\G)(\k)$, on r{\'e}sout un syst{\`e}me lin{\'e}aire (avec second
membre) {\`a} coefficients rationnels. Le syst{\`e}me homog{\`e}ne associ{\'e} est
exactement celui qui d{\'e}finit la composante homog{\`e}ne de poids $n+1$ de 
$\dmrd(\G)(\k)$. La diff{\'e}rence de deux solutions est donc un {\'e}l{\'e}ment
de $\dmrd(\G)(\k)$. D'autre part, modulo des termes de poids au moins
$n+2$, si $\psi$ est homog{\`e}ne de poids $n+1$, l'op{\'e}rateur
$\exp(s_\psi)$ agit sur une s{\'e}rie de terme constant 1 par addition de
$\psi$. On en tire par une r{\'e}currence facile l'{\'e}nonc{\'e} suivant~:
\begin{prop}\label{prop_trans}
Si $\DMRD_\l(\G)(\k)$ est non-vide, l'action de $\Exp(\dmrd_0(\G)(\k))$ sur
chaque troncation $\tr{\DMRD_\l}{n}(\G)(\k)$ et sur $\DMRD_\l(\G)(\k)$ est
transitive.  
\end{prop}

Pour tout $\mu\in\k$, soit $h_\mu$ le morphisme de 
$\k$-alg{\`e}bres topologiques multipliant chaque lettre de $X_G$ par
$\mu$. On voit facilement que $\DMRD(\G)(\k)$ est stable par $h_\mu$. 
Lorsque le cardinal de $\G$ vaut au moins
3, l'{\'e}l{\'e}ment exceptionnel $\alpha_\G$ est de poids 1. On a donc
$h_\mu(\DMRD_\l(\G)(\k))\subset\DMRD_{\l\mu}(\G)(\k)$. Si $\G=\triv$
ou $\{\pm1\}$, on a
$h_\mu(\DMRD_\l(\G)(\k))\subset\DMRD_{\l\mu^2}(\G)(\k)$, car
$\alpha_G$ est de poids 2.  
Comme $\LC_\ssh(\G)$ appartient {\`a} $\DMRD(\CM)$, en posant, suivant le
cas, $\mu=(\LC_\ssh(\G)|\alpha_G)^{-1}$ ou
$\mu=(\LC_\ssh(\G)|\alpha_G)^{-2}$, on obtient par $h_\mu$ un {\'e}l{\'e}ment de
$\DMRD_1(\G)(\CM)$, qu'on notera $\ov\LC_\ssh$.  
\begin{prop}
   Il existe un {\'e}l{\'e}ment $\Psi$ de $\DMRD_1(\G)(\QM)$.
\end{prop}
On r{\'e}sout pour cela par r{\'e}currence les {\'e}quations d{\'e}finissant
$\DMRD_1(\G)$. Si $\Phi$ est un {\'e}l{\'e}ment de $\tr{\DMRD_1}{n}(\G)(\QM)$,
il existe $\psi\in\dmrd_0(\G)(\CM)$ tel que
$\Theta\ \ass\ \exp(s_\psi)(\ov\LC_\ssh(\G))$ soit {\'e}gal {\`a} $\Phi$
modulo des termes de poids au moins $n+1$. Or
$\Phi+\Theta_{n+1}=\tr{\pi}{n+1}(\Theta)$ 
appartient {\`a} 
$\tr{\DMRD_1}{n+1}(\CM)$, en notant $\Theta_{n+1}$ le terme de poids
$n+1$ de $\Psi$. Le syst{\`e}me lin{\'e}aire satisfait par les coefficients
d'un {\'e}l{\'e}ment $\Phi_{n+1}$ homog{\`e}ne de poids $n+1$ de $\assl{\QM}{X_\G}$ pour
que $\Phi+\Phi_{n+1}$ appartienne {\`a} $\tr{\DMRD_1}{n+1}(\G)(\QM)$ est
enti{\`e}rement rationnel et admet une solution complexe. Il admet donc
une solution rationnelle. 

Si le cardinal de $\G$ vaut au moins 3, on en d{\'e}duit 
que $\DMRD_\l(\G)(\k)$ est non-vide pour tous $\l$ et $\k$~: il
contient $h_\l(\Psi)$. Dans
l'autre cas, on montre d'abord l'existence d'un {\'e}l{\'e}ment 
$\Psi_{\text{pair}}$ de $\DMRD_1(\G)(\QM)$ dont tous les termes de
poids impairs sont nuls ~: avec les notations ci-dessus, $0$ est 
solution du syst{\`e}me
lin{\'e}aire portant sur $\Phi_{n+1}$ si $n+1$ est impair tandis que $\Phi$ est
pair. On utilise alors $h_\l^{\text{pair}}$, 
l'op{\'e}rateur qui multiplie les mots de poids $2n$ par $\l^n$. La
proposition \ref{prop_trans} permet alors de conclure. 
\section{Cons{\'e}quences et conclusion}\label{sec_conseq}
Comme corollaire, l'application de $\k\times\dmrd_0(\G)(\k)$ dans
$\DMRD(\G)(\k)$ donn{\'e}e par $(\l, \psi)\mapsto \Exp(\psi)\pmt
h_\l(\Psi)$ (ou, suivant le cas, $\Exp(\psi)\pmt
h_\l^{\text{pair}}(\Psi_{\text{pair}})$) est un isomorphisme
de sch{\'e}mas de $\AM^1\times \dmrd_0(\G)$ sur $\DMRD(\G)$. 
Le foncteur
$\k\mapsto\dmrd_0(\G)(\k)=\dmrd_0(\G)(\QM)\whot\k$ est repr{\'e}sentable par
l'alg{\`e}bre sym{\'e}trique du dual gradu{\'e} de $\dmrd(\G)(\QM)\cap\assl{\QM}{X_\G}$. 
Cela montre que l'alg{\`e}bre formelle $\DC\MC\RC\DC(\G)$ d{\'e}finie par les
relations DMRD est une 
alg{\`e}bre de polyn{\^o}mes car $\DMRD(\G)$ est son spectre et donne une
description de ses g{\'e}n{\'e}rateurs libres.

La libert{\'e} de l'alg{\`e}bre $\DC\MC\RC\DC(\triv)$ constitue une
partie de l'{\'e}nonc{\'e} du th{\'e}or{\`e}me d'{\'E}calle. Il donne une description des
g{\'e}n{\'e}rateurs en fonction de son alg{\`e}bre de Lie des polyn{\^o}mes
bialternaux, elle-m{\^e}me sous-alg{\`e}bre de Lie de ARI. Il a tr{\`e}s r{\'e}cemment
{\'e}tendu ces r{\'e}sultats au cas $\G=\{\pm1\}$. Une fois {\'e}tablie la
correspondance entre les s{\'e}ries g{\'e}n{\'e}ratrices {\em commutatives} que
Goncharov et lui utilisent et les s{\'e}ries g{\'e}n{\'e}ratrices non-commutatives
de cette note, on constate que $\mt(\G)$ est isomorphe {\`a} la
sous-alg{\`e}bre d'ARI form{\'e}e des \og{}moules entiers\fg{}, \ie{}
correspondant {\`a} des s{\'e}ries enti{\`e}res au voisinage de 0 et dont les
variables $u$ sont dans $\G$, cas
particulier dont les formules sont donn{\'e}es chez Goncharov 
(\cite{Gonch98}). Bien que d{\'e}velopp{\'e}s ind{\'e}pendamment, les arguments
utilis{\'e}s par {\'E}calle pour conclure semblent 
fortement recouper les n{\^o}tres (lin{\'e}arisation et utilisation de la stabilit{\'e}
pour le crochet de $\mt(\G)$ des alg{\`e}bres de Lie concern{\'e}es). 

Suivant la variante de Drinfel'd de la conjecture de Deligne
(\cite{DrinQTQH}, p.860, questions), l'alg{\`e}bre de Lie $\grt_1$ du groupe $\Ass_0=\GRT_1$
(lequel agit sur les 
associateurs par translation {\`a} gauche au sein de $\MT(\triv)$) est
une alg{\`e}bre de Lie libre, avec un g{\'e}n{\'e}rateur et un seul en chaque
poids impair except{\'e} 1. Drinfel'd exhibe, apr{\`e}s Ihara, un syst{\`e}me
d'irr{\'e}ductibles de $\grt_1(\CM)$, avec les bonnes conditions de poids,
en les lisant dans $\Phi_\KZ$ (qui est {\'e}gal {\`a} $\LC_\ssh(\triv)$). Plus 
pr{\'e}cis{\'e}ment, ce sont les composantes homog{\`e}nes de l'{\'e}l{\'e}ment $\psi$ de
$\grt_1(\CM)$ tel que
$h_{-1}(\Phi_\KZ)=\Exp(\psi)\pmt\Phi_\KZ$. D'apr{\`e}s le th{\'e}or{\`e}me
\ref{theo_action}, $\psi$ appartient {\'e}galement {\`a} $\dmrd(\triv)(\CM)$,
ainsi donc que les irr{\'e}ductibles de Drinfel'd. Si la conjecture est
vraie, on doit donc avoir $\GRT_1\subset\DMRD(\triv)$. D'autre part,
la conjecture de Zagier pr{\'e}voit  la dimension de l'espace vectoriel
des polyz{\^e}tas de poids $n$. On peut formuler, pour s'affranchir des
probl{\`e}mes de transcendance, une variante formelle, \ie{} portant sur
la dimension des composantes homog{\`e}nes de
$\DC\MC\RC\DC(\triv)$. Si cette conjecture est {\'e}galement vraie,
par {\'e}galit{\'e} des dimensions, l'inclusion ci-dessus est une
{\'e}galit{\'e}. Ceci motive la conjecture \ref{conj} {\'e}nonc{\'e}e {\`a} la section
\ref{sec_tors}. 
\section{Tables de dimensions}
Les dimensions ci-dessous ont {\'e}t{\'e} calcul{\'e}es par ordinateur. {\`A} titre
indicatif, on donne {\'e}galement les dimensions de $\dmr(\Gamma)$,
l'espace tangent {\`a} $\DMR(\Gamma)$ au voisinage de 1. Par $\mu_N$, on
entend le groupe des racines $N$\emes{} 
de l'unit{\'e} dans $\CM$. Les symboles $\dag$ signalent la pr{\'e}sence d'un 
{\'e}l{\'e}ment irr{\'e}gulier. 

\vspace{1em}
\begin{center}
Composante homog{\`e}ne de poids $p$ de $\dmr(\triv)=\dmrd(\triv)$ 
\input{dimdmr1.tex}
\medskip

\begin{tabular}{p{5cm}p{2eM}p{5cm}}
Composante homog{\`e}ne de poids $p$ de $\dmrd(\mu_N)$ && 
Composante homog{\`e}ne de poids $p$ de $\dmr(\mu_N)$ \\
\input{dimdmrd.tex}&&
\input{dimdmr.tex}
\end{tabular}

\end{center}

\nocite{IhICM, HoffAlg, Demaz, IhIsrael, Broad}

\bibliographystyle{joeplain}
\bibliography{../../bibtex/these,../moi}

\end{document}

%% file: dimdmr1.tex

\begin{tabular}{|c|*{11}{c}|}
\hline $N\backslash p$&1&2&3&4&5&6&7&8&9&10&11\\
\hline 1 &0 &1$^\dag$ &1 &0 &1 &0 &1 &1 &1 &1 &2\\
\hline\end{tabular}

%% file: dimdmrd.tex

\begin{tabular}{|c|*{7}{c}|}
\hline $N\backslash p$&1&2&3&4&5&6&7\\
\hline 2 &1 &1$^\dag$ &1 &1 &2 &2 &4\\
3 &2$^\dag$ &1 &2 &3 &6 & &\\
4 &2$^\dag$ &1 &3 &6 & & &\\
5 &3$^\dag$ &2 &6 &13 & & &\\
6 &3$^\dag$ &2 &7 & & & &\\
7 &4$^\dag$ &4 &13 & & & &\\
8 &3$^\dag$ &4 &15 & & & &\\
9 &5$^\dag$ &7 &23 & & & &\\
10 &4$^\dag$ &6 &26 & & & &\\
11 &6$^\dag$ &10 & & & & &\\
\hline\end{tabular}

%% file: dimdmr.tex

\begin{tabular}{|c|*{7}{c}|}
\hline $N\backslash p$&1&2&3&4&5&6&7\\
\hline 2 &1 &1$^\dag$ &1 &1 &2 &2 &4\\
3 &2$^\dag$ &1 &2 &3 &6 & &\\
4 &3$^\dag$ &1 &3 &7 & & &\\
5 &4$^\dag$ &2 &6 &13 & & &\\
6 &5$^\dag$ &3 &8 & & & &\\
7 &6$^\dag$ &4 &13 & & & &\\
8 &7$^\dag$ &5 &17 & & & &\\
9 &8$^\dag$ &7 &23 & & & &\\
10 &9$^\dag$ &8 &31 & & & &\\
11 &10$^\dag$ &10 & & & & &\\
\hline\end{tabular}